\documentclass[preprint,11pt]{elsarticle}

\usepackage{amssymb}

\usepackage{natbib}
\usepackage{graphicx}
\usepackage{amsmath,amsfonts}
\usepackage{arydshln}
\usepackage{xcolor}
\usepackage{ragged2e}

\journal{Mechanical Systems and Signal Processing}

\begin{document}
\title{Using scientific machine learning for experimental bifurcation analysis of dynamic systems}
\author{Sandor Beregi}
\author{David A. W. Barton}
\author{Djamel Rezgui}
\author{Simon Neild}

\address{Faculty of Engineering, University of Bristol, United Kingdom}

\begin{abstract}
    Augmenting mechanistic ordinary differential equation (ODE) models with machine-learnable structures is an novel approach to create highly accurate, low-dimensional models of engineering systems incorporating both expert knowledge and reality through measurement data. Our exploratory study focuses on training universal differential equation (UDE) models for physical nonlinear dynamical systems with limit cycles: an aerofoil undergoing flutter oscillations and an electrodynamic nonlinear oscillator. We consider examples where training data is generated by numerical simulations, whereas we also employ the proposed modelling concept to physical experiments allowing us to investigate problems with a wide range of complexity. To collect the training data, the method of control-based continuation is used as it captures not just the stable but also the unstable limit cycles of the observed system. This feature makes it possible to extract more information about the observed system than the open-loop approach {(surveying the steady state response by parameter sweeps without using control)} would allow. We use both neural networks and Gaussian processes as universal approximators alongside the mechanistic models to give a critical assessment of the accuracy and robustness of the UDE modelling approach. We also highlight the potential issues one may run into during the training procedure indicating the limits of the current modelling framework.
    %\keywords{machine learning \and nonlinear dynamics \and universal differential equations \and bifurcation analysis \and aeroelastic flutter}
\end{abstract}

\maketitle

\section{Introduction}
\label{intro}

Data-driven modelling of real-life systems is an increasingly popular approach both in the scientific community and in industry. With increasing data and computational capacity becoming available, there is a growing number of areas where machine learning techniques are used to incorporate measurement data in the modelling procedure \cite{SciMLWSrep}. This also resonates with the industrial demand for highly precise, adaptive, real-time models of physical structures. These so-called digital twins \cite{WaggDigiTwin} can then be used for example in design, decision-making, or failure diagnostics.

While machine learning techniques have been successfully used in many fields such as image processing \cite{Krizhevsky2017}, speech recognition \cite{Graves2013} or computational biology \cite{CompBio}, the field of scientific machine learning (SciML) is relatively new \cite{SciMLWSrep}. The feature of SciML that differentiates it from most other machine learning applications is that machine learning techniques are employed alongside existing mathematical models; i.e., instead of treating the modelled system as an unknown `black box'-type system, additional information such as symmetries or known physical laws are incorporated in the hybrid machine-learned models \cite{Raissi2019}. {This modelling strategy falls into the scale of `grey-box'-type models where mechanistic, `white-box' modelling is augmented rather than replaced by data-driven techniques \cite{Worden2018}. The general idea behind this concept is to exploit the advantages of both (black- and white-box-type) modelling approaches: achieving high fidelity to measurement data, while relying on the physics-based core of the model to assist in extrapolation to parameter-domains not covered by the training data -- an area often posing a challenge for pure data-driven models -- and reducing the scale of the required training dataset \cite{Pitchforth2021} .}

{In \cite{Kennedy2001}, Kennedy and O'Hagan introduce this concept for a generic computer code, taking a Bayesian approach to calibrate the model parameters and correcting its error to measurement data. Worden et al. \cite{Worden2018} consider this approach applied to nonlinear systems, whereas Pitchforth et al. \cite{Pitchforth2021} present a study on training grey box models capable of predicting wave-load on offshore structures. }

{Our study considers a similar modelling concept applied specifically to replicate the behaviour of observed nonlinear systems from their measured bifurcation diagrams and the associated periodic solutions. The complexity of this task introduces several challenges, the above mentioned works do not specifically account for. One such challenge is that a bifurcation diagram is not just a mere collection of solutions of the respective system, but it also includes relevant information about the stability of the steady states it represents. Thus, by topological extrapolation it allows for the qualitative prediction of the transient behaviour of the system. This means that e.g. an accurate prediction of a limit cycle does not automatically imply a good representation of the observed system since the stability of the solution may be wrongly captured. Therefore, the fitness of these models may not be solely assessed based on the error to the measurements. This also highlights the advantage of having a qualitatively correct physics-based core in the model since these features of the dynamic system would be very challenging to capture in a pure data-based manner.}

Using universal differential equations (UDEs), that is, differential equations with embedded machine-learnt structures is a recent approach within SciML. UDEs are essentially hybrid differential equations with machine learnable structures embedded in the right-hand-side. Since differential equations are commonly used to model engineering problems, UDEs provide a convenient way to incorporate expert knowledge and physical insight about the observed system into machine-learnt models. This modelling concept is similar to neural differential equations (NDEs), {which can be interpreted as a subset of UDEs} where the derivative of the state variables is modelled by a neural network \cite{Chen2018}. Nevertheless, making `smaller' corrections to qualitatively correct physics-based models, rather than relying purely on machine-learnable structures to capture the dynamics of the observed system, is expected to require a smaller computational effort by enabling to use smaller-scale machine-learnt structures and smaller data-sets.

With physics-based models on the other hand, it is often not feasible to consider and capture all the experimental features within. Thus, mechanistic models in almost every case have some degree of error compared to measurement data, particularly when nonlinear behaviour {and a complex structure of steady state solutions are} present. Thus, the accuracy of these models can be further improved by machine-learnt structures trained on measurement data.

Several studies have used UDE or NDE models for physical structures \cite{RackauckasHelicopter}, pandemic- \cite{DandekarCovid,Nunez} and climate modelling \cite{Ramadhan2020}. It is worth highlighting the work of Rackauckas et al. \cite{DiffEqFlux} to create an environment in the Julia programming language which allows for construction and training of such models. Most of these studies focus on identifying a well-fitting model with a constant set of system parameters. Nevertheless, there are several applications where it is also interesting to investigate the system's response to varying parameters. For instance, one may be interested in the effect of parameter-uncertainty on the system, or its behaviour under changing external conditions such as temperature, or wear during its life-cycle.

Our study focuses on another such area, nonlinear dynamical systems and the application of the UDE models in this context. Nonlinear dynamical systems are commonly characterised by their bifurcation diagrams indicating the steady-state system behaviour as selected system parameters (referred to as bifurcation parameters) change. Thus, these are typically considered as varying-parameter problems \cite{Kuznetsov}.

While obtaining bifurcation diagrams for numerical models is standard nowadays and aided by numerous software packages, doing the same for physical experiments is much more recent. The main issue in experiments is that by varying one of the system parameters and observing the response in an `open-loop' way, one can only identify the stable solutions, whereas the unstable ones remain hidden. The technique of control-based continuation (CBC), introduced in \cite{SieberCBC} is capable of tracing both stable and unstable limit cycles in physical experiments. The method has been adapted to deal with random perturbations/process-noise in the experiment \cite{SCHILDER2015,Renson2019}, it has been also extended to provide information on stability (without turning the controller off) \cite{Barton2017}, and it has been made capable of tracing self-excited vibrations in autonomous systems \cite{Lee2020}. Having the additional information about the unstable solutions has significant merits in modelling and parameter identification problems. As shown in \cite{cbc_paper}, by identifying the parameters of a dynamical system with polynomial nonlinearity, CBC, when compared to the {open-loop approach}, offers better robustness against process noise in revealing the finer details of physical experiments. This result motivated the use of CBC in this study to collect the training data for machine-learnable models.

In context of nonlinear systems, there is often an interest to find the so-called normal forms: {simplified} models with the bifurcation diagram characteristic to the observed phenomenon. In their study, Lee et al. \cite{Lee2020} used machine learning techniques to identify the transformation from physical to normal coordinates of a Hopf bifurcation for aeroelastic flutter based on measurement data. The aim of our investigation is similar in the sense that we reconstruct measured bifurcation diagrams of nonlinear systems: an aerofoil undergoing flutter and an electrodynamical oscillator. However, instead of identifying the reduced order normal form, we do not restrict the model to a minimal order, but instead, use as many of the available states from the measurement in the UDE model as possible while retaining the physical coordinates.

Our study is exploratory in the sense that our main goal is to assess the potential of the mechanistic/machine-learnt hybrid models for nonlinear systems. Therefore, we choose the approach of using relatively simple machine-learnt models and training framework and reporting on the issues that occurred, rather than investing effort to tailor the algorithms to tackle these challenges.

The rest of our paper is organised as follows. In Section 2, we introduce the modelling concept using UDE models comprising a physics-based model with the measured state variables of the system and a machine-learnable structure trained to compensate the error of the pure mechanistic model. In Section 3, two examples are presented, an aerofoil that exhibits aeroelastic flutter, and nonlinear electrodynamic oscillator. In Sections 4 and 5, training of UDE models is performed using data from numerical models and physical experiments, respectively. For our two practical examples, we draw the conclusions of our analysis in Section 6.

\section{Modelling concept}

\begin{figure*}
  \includegraphics[width=\textwidth]{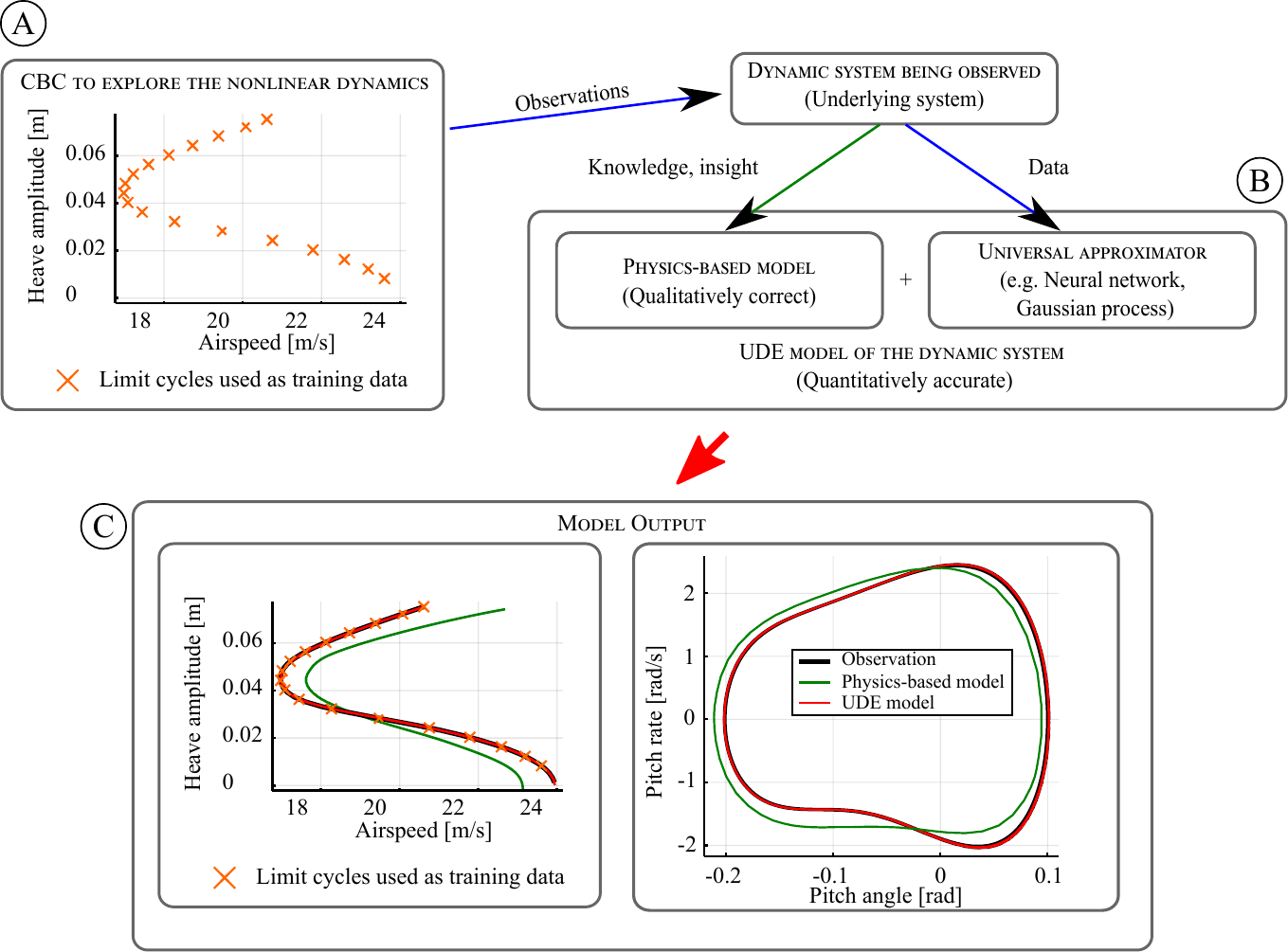}
\caption{Panel (A): Collecting measurement data using control-based continuation. Panel (B) The concept of physics-based/machine-learnable hybrid models for nonlinear dynamical systems. Panel (C): Model output of the UDE model with improved prediction of the bifurcation diagram and the trajectories.}
\label{fig:concept}       % Give a unique label
\end{figure*}

Our study focuses on training UDE models for nonlinear systems exhibiting limit cycles. It is common practice to represent these in bifurcation diagrams where branches of limit cycles are shown while one (or more) system parameter(s) -- the bifurcation parameter(s) -- is/are varied \cite{Kuznetsov}. These diagrams are fundamental in characterising nonlinear systems since they not only show the steady state solutions and their stability, but, in addition, with some topological extrapolation, give information about the transient behaviour. These characteristics motivate us to use measured limit cycles in the numerical simulations and physical experiments both as training data and as benchmark to evaluate the predictions of the fitted UDE models. The modelling concept is summarised in Fig.~\ref{fig:concept}.

{We consider a nonlinear experiment as the observed system where the steady state behaviour, represented by bifurcation diagrams, can be measured using control-based continuation. Our study focuses specifically on applications with bistable behaviour; nevertheless, the modelling approach may be useful to find a representative model of a wider group of dynamical systems. The general idea is that based on knowledge and insight into the behaviour of the observed system, one can derive a qualitatively correct physics-based model using the available (measured) state variables. Then, measured bifurcation diagrams, i.e., the time-profiles of the limit cycles associated with the measured points in the bifurcation diagrams are used with machine-learnt structures to account for the unconsidered effects affecting the observed system. The expected output of the UDE model is a both qualitatively and quantitatively accurate numerical representation of the observed system capable of replicating the measured bifurcation diagram and the corresponding trajectories.}

\subsection{The universal differential equation model}

Let us consider an experimental representation of a dynamical system. In real-life scenarios, one can only acquire data from the observable part of the underlying dynamical system and its state variables. Therefore, we divide the \emph{observed system} and its state variables into a \emph{measured} and an {\emph{unobserved}} part. We refer to the vectors of the corresponding measured and {unobserved} state variables as $\mathbf{x}$ and $\mathbf{u}$, respectively, with $\mathbf{x}(t) \in \mathcal{R}^m$ and $\mathbf{u}(t) \in \mathcal{R}^{m_{\rm h}}$, {where $m$ is the dimension of the measured part of the model whereas $m_{\rm h}$ accounts for the number of state variables not captured in the experiment. The model parameters are given in the vector $\mathbf{p}$. }. To take into account the effect of process noise, the underlying system can be given as a stochastic differential equation in the form
\begin{equation}\label{eq:modeleq}
\left( \begin{matrix}
  \mathrm{d}\mathbf{x} \\
  \hdashline
  \mathrm{d}\mathbf{u} \\
\end{matrix} \right) = \left( \begin{matrix}
  \mathbf{f}(\mathbf{x}, \mathbf{u}; \mathbf{p}) \\
  \hdashline
  \mathbf{g}(\mathbf{x}, \mathbf{u}; \mathbf{p}) \\
\end{matrix} \right) \mathrm{d}t + \left( \begin{matrix}
  \boldsymbol{\Gamma} (\mathbf{x}, \mathbf{u}; \mathbf{p}) \\
  \hdashline
  \boldsymbol{\Phi} (\mathbf{x}, \mathbf{u}; \mathbf{p}) \\
\end{matrix} \right) \mathrm{d}W, \quad \begin{matrix}
  \rm{measured} \\
  \hdashline
  \rm{unobserved}
\end{matrix}
\end{equation}
where the functions $\mathbf{f}$ and $\mathbf{g}$ define the deterministic parts of the system. {We assume that this deterministic core is sufficiently smooth in the observed domain(s) of the parameters and the state variables. In an experiment, the assessment of smoothness is often less straightforward. In principle, many phenomena which are usually present in a physical test to a certain degree, e.g. friction or freeplay in joints have non-smooth mechanical models. Nevertheless we assume that the effect of non-smoothness is sufficiently small for the experiments to be approximated by smooth models. }

The effect of process noise is considered by the Wiener process $W$ with the coefficient functions $\boldsymbol{\Gamma}$ and $\boldsymbol{\Phi}$ corresponding to the measured and {unobserved} parts of the system, respectively. This assumes a single source of noise in the observed system, which covers all examples considered in this study. Nevertheless, the concept could be generalised to accommodate multiple sources of noise. {Note that the role of the Wiener process in the model is similar to the assumed unconsidered state variables, accounting for an unknown part of the observed system dynamics. This may seem redundant but the two approaches have slightly different purpose. While the unconsidered part accounts for the limitations of the number of measured variables one can obtain in an experiment, fluctuations caused by an unknown or difficult to model phenomenon may appear in the measured signals; thus, justifying the use of a stochastic model. }

Equation \eqref{eq:modeleq} simplifies to an ordinary differential equation (ODE) with $\Gamma = \mathbf{0}$ and $\boldsymbol{\Phi} = \mathbf{0}$ which we refer to as the deterministic case.

We approximate the measured part of the underlying system with a universal differential equation \cite{DiffEqFlux} where the \emph{physics-based part} $\tilde{\mathbf{f}}$ is augmented with a so-called \emph{universal approximator}: a machine-learnable regression function $\mathbf{U}$ with parameters $\mathbf{q}$
\begin{equation}\label{eq:ude_app}
    \mathrm{d} \mathbf{x} \approx \mathrm{d} \tilde{\mathbf{x}} := \left( \tilde{\mathbf{f}} (\tilde{\mathbf{x}}; \mathbf{p})  +  {\mathbf{U}} (\tilde{\mathbf{x}}; \mathbf{p}; \mathbf{q}) \right)\rm{d} t +  \tilde{\boldsymbol{\Gamma}} (\tilde{\mathbf{x}}; \mathbf{p}) \rm{d} W.
\end{equation}
Accordingly, it is assumed that the physics-based part and the undescribed dynamics are separable. It is worth mentioning that the modelling framework would allow to consider a universal approximator embedded to the physics-based part $\tilde{\mathbf{f}}$ as $\tilde{\mathbf{f}} \left( \tilde{\mathbf{x}}, \mathbf{U} (\tilde{\mathbf{x}}; \mathbf{p}; \mathbf{q}); \mathbf{p} \right)$ which can be useful in scenarios when one has some -- albeit limited -- insight into the dynamics that are missing from the mechanistic model.

We assume that the physics-based part of the model already provides qualitatively correct results by itself; hence, the main purpose of the correction with the universal approximator is to compensate the quantitative difference between the theoretical results and measurements. While it would be possible to approximate the entire right-hand side of Eq. \eqref{eq:ude_app} with universal approximators, in doing this we would fail to exploit any prior knowledge about the modelled system. The main benefit of this hybrid approach is that a smaller-scale universal approximator and/or fewer measurements may be satisfactory to obtain an accurate numerical model than with a pure machine-learnable model. This might be especially useful if data acquisition from the measurement is expensive.

In our study, two of the most commonly used universal approximators are considered: \emph{neural network}s and \emph{Gaussian process}es.

\subsection{Model training}

{In this subsection below, we provide an introductory description on neural networks and Gaussian processes and the associated model training procedure. Note that, in this study, we use the standard types of the machine-learnt structures, relatively shallow, feed-forward neural networks and Gaussian processes with an SE-ARD kernel. This is intentional since the aim was to test the capability of a generic version of these structures to augment the physics-based models. Nevertheless, we point to more special techniques that might be used to address the challenges encountered during the training procedure. The notations in the presented formulae were derived from \cite{Higham2019} for neural networks and \cite{GaussianProcessesJL} for Gaussian Processes. For further details on the employed machine-learnt structures, we invite the reader to follow the references provided.}

\subsubsection{Neural networks}

When using a neural network as the universal approximator, we employ the DiffEqFlux.jl \cite{DiffEqFlux} package in Julia to train the neural network in the UDE model. This package combines numerical simulations with machine learning, as the training procedure is based on minimising the error between the measured and predicted trajectories.

  Let us consider the $j$-th approximation of the measured part of the observed system in the training procedure, using the parameters $\mathbf{q}_j$ in the universal approximator
\begin{equation}\label{eq:ude_app_j}
    \mathrm{d} \tilde{\mathbf{x}}_j := \left( \tilde{\mathbf{f}} (\tilde{\mathbf{x}}; \mathbf{p})  +  {\mathbf{U}} (\tilde{\mathbf{x}}; \mathbf{p}; \mathbf{q}_j) \right)\rm{d} t +  \tilde{\boldsymbol{\Gamma}} (\tilde{\mathbf{x}}; \mathbf{p}) \rm{d} W.
\end{equation}
We refer to the predicted trajectories corresponding to the initial condition ${\tilde{\mathbf{x}}_{0}}$ and parameters $\mathbf{p}$ and $\mathbf{q}_j$ as $\tilde{\mathbf{x}}_j(t)$. Assuming we have measurement data collected from the observed system \eqref{eq:modeleq} at the time instants $t_i$, $i=1,...,N$, the measured trajectories are represented by $\mathbf{x}(t_i) =: \mathbf{x}_i$. The objective function is formulated using the Euclidean norm of the sum of errors between the measured and predicted states
\begin{equation}
  \varepsilon^j = \sum_{j = 1}^N \left( \tilde{\mathbf{x}}_i^j - \mathbf{x}_{i} \right)^2 .
\end{equation}

The model parameters can be divided into two groups: the parameters $\mathbf{p}$ of the physics-based part and the parameters $\mathbf{q}$ corresponding to the neural network. We assume that the parameters of the physics-based model are already optimised and fitted to match the measurement data as well as possible. This prevents the neural network from approximating existing known terms. Therefore, our analysis focuses on the optimisation of the neural network parameters only.

In our study, we use feed-forward neural networks \cite{LeCun2015,Shreshta2019NN,Higham2019} consisting of nodes with a sigmoid activation function. The nodes are organised in layers which are referred to either as \emph{input}, \emph{hidden} or \emph{output} layers according to their position in the network. In a neural network with $L$ layers (the 1$^{\rm st}$ layer is referred to the as input-, the $L$-th as the output layer) the contribution of the $k^{\rm th}$ layer can be given as
\begin{equation}
\mathbf{y}_k = \mathbf{h}_{\rm act}(\mathbf{Z}_k) = \mathbf{h}_k(\mathbf{y}_{k-1}),
\end{equation}
where $\mathbf{h}_{\rm act}$ is the activation function with $\mathbf{y}_k \in \mathbf{R}^{N_k}$, $\mathbf{y}_{k-1} \in \mathbf{R}^{N_{k-1}}$ while
\begin{equation}
  \mathbf{Z}_k = \mathbf{W}_k \mathbf{y}_{k-1} + \mathbf{b}_k,
\end{equation}
with biases $\mathbf{b}_k \in \mathbf{R}^{N_k}$ and weights $\mathbf{W}_k \in \mathbf{R}^{N_k \times N_k}$.
Thus, the effect of the whole neural network can be given by the nested function
\begin{equation}
\mathbf{y}_L = \mathbf{h}_L \circ \mathbf{h}_{L-1} \circ ... \circ \mathbf{h}_1 (\mathbf{y}_0),
\end{equation}
where $L$ is the number of layers in the network.
Thus, the aim of the optimisation procedure is to find the elements of the weight matrices $\mathbf{W}_k$ and biases $\mathbf{b}_k$ of the nodes in the neural network corresponding to the global minimum of the objective function $\varepsilon$.

The training of the neural network is carried out in multiple-steps using the GalacticOptim.jl \cite{GalacticOptim} package in Julia. The ADAM algorithm in the Flux.jl Julia package \cite{Flux} is used first in our examples to find an approximate minimum of the objective function. Then, the identified minimum is refined by standard gradient descent.

%While this algorithm is fairly effective in the initial steps, it is prone to identifying local minima as the optimum. To avoid this, the application of regularisation techniques may be necessary.

\subsubsection{Gaussian processes}

An alternative approach to neural networks is to use a Gaussian process as the universal approximator \cite{GaussianBook}. In the literature, several studies investigate the relationship between neural networks and Gaussian processes, finding that certain neural networks have an equivalent Gaussian process \cite{Lee2017,Rudner2017}.

Since Gaussian process regression directly uses the training dataset to make predictions, the regression function is optimised using the right-hand side of \eqref{eq:ude_app}. Let us assume we have access to the derivatives $\rm{d}\mathbf{x}_t$ in the time instants $t_i$: $\rm{d} \mathbf{x}_i$. Accordingly, the error in the right-hand side is expressed as
\begin{equation}
\varepsilon_{\rm acc}^j = \sum_{i=1}^N \left( {\rm d} \mathbf{x}_i - {\rm d} \tilde{\mathbf{x}}_i^j \right)
\end{equation}

The training, however, is not solely based on the error between the prediction and the measurement data. A significant difference between the two approaches is that while neural networks are deterministic, Gaussian processes have a stochastic underlying function characterised by the mean $m_{\boldsymbol{\theta}}(\mathbf{x})$ and the covariance $k_{\boldsymbol{\theta}}(\mathbf{x}, \mathbf{x}')$ determined by the kernel function and the inputs $\mathbf{X} = (\mathbf{x}_1, ... , \mathbf{x}_n)$. In our study, we used a squared exponential kernel with automatic relevance prediction (SE-ARD). {The advantage of using an SE-ARD kernel compared to a standard squared exponential kernel is its capability of handling training data where the different parameters vary on a significantly different length-scale. The additional parameters in the SE-ARD kernel allow for} considering different length-scales for the input variables. The SE-ARD kernel function is given by
\begin{equation}
k_{\boldsymbol{\theta}}(\mathbf{x}, \mathbf{x}') = \sigma^2 \rm{exp} ((\boldsymbol{\lambda} \odot (\mathbf{x}-\mathbf{x}'))^2),
\end{equation}
with $\boldsymbol{\lambda} = \left( \lambda_1,...,\lambda_m  \right)$, where $\lambda_i = 1/(2 \ell_i)$, while $\odot$ denotes the Hadamard product. Thus, the vector $\boldsymbol{\theta}$ of the hyperparameters contains the length scales $\ell_i$ of the kernel function and the standard deviation $\sigma$. Consequently, the number of parameters to be optimised in a Gaussian process is much lower than in case of a neural network which is a significant advantage of using the Gaussian process.

Let us assume we have a dataset $\mathfrak{D} = (\mathbf{X},\mathbf{Y})$ with observations $\mathbf{Y} = (\mathbf{y}_1, ... , \mathbf{y}_n)$ and covariates
 \noindent $\mathbf{X} = (\mathbf{x}_1, ... , \mathbf{x}_n )$. Then, the Gaussian Process fit is assessed by the marginal likelihood $p(\mathfrak{D} \vert \boldsymbol{\theta}, \sigma)$, i.e., we would like to find $\theta$ corresponding to the maximal likelihood. In optimisation algorithms, the negative logarithm \\
 \noindent $ - {\rm log} \left( p(\mathfrak{D} \vert \boldsymbol{\theta}, \sigma) \right)$ is used as the objective function.

The optimisation is realised in Julia using the GaussianProcesses.jl package \cite{GaussianProcessesJL}. Due to the small number of parameters, optimising the parameters of a Gaussian process takes a fraction of the time needed to train a neural network. However, since the Gaussian process needs to use the training data for predictions, it does not scale well if the training dataset is large. Moreover, one may also experience a sensitivity to the initial guess during the optimisation which might make it necessary to find an appropriate initial guess by iteration. Therefore, Gaussian processes are most effective with relatively simple problems in terms of size and dimension.

\section{Practical examples}

We introduce two practical use-cases for the UDE models: aeroelastic flutter and a forced electrodynamical oscillator; both systems featuring a nonlinear restoring force.

\subsection{Aerofoil prone to aeroelastic flutter}\label{sec:Flutter_num}
\begin{figure}
\begin{center}
  \includegraphics[width=0.9\textwidth]{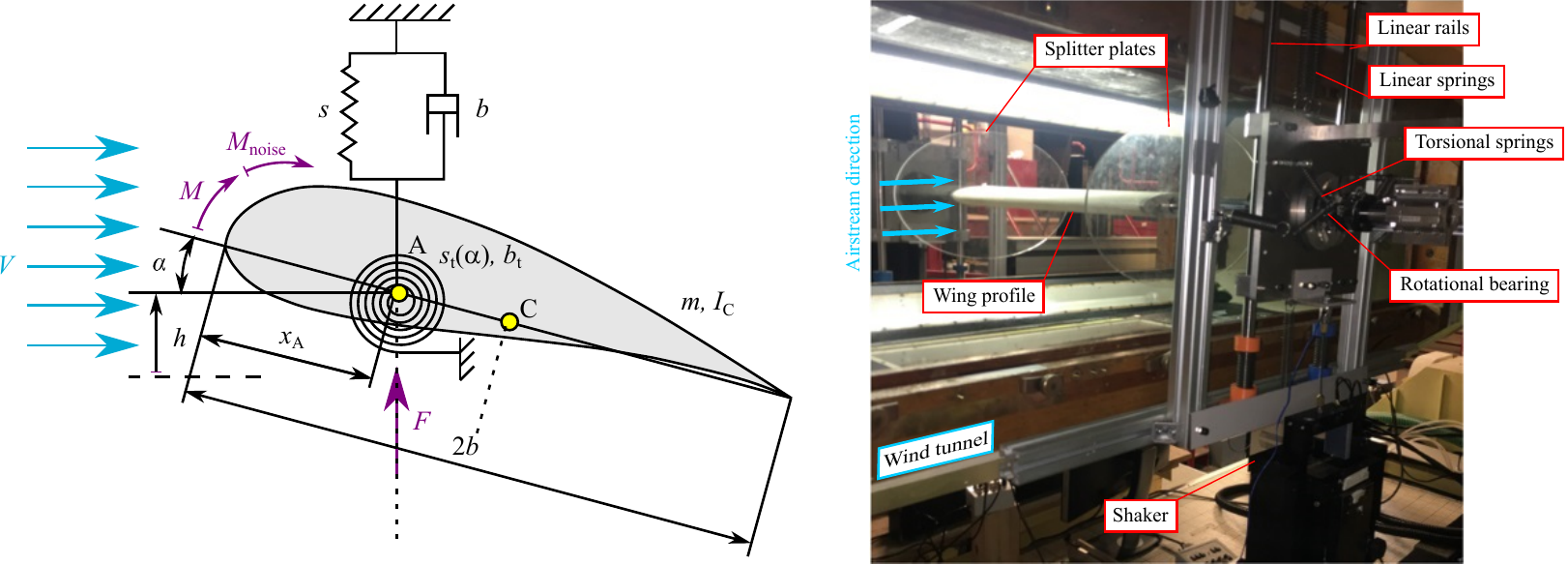}
\end{center}
\caption{Left: Model of an elastically supported aerofoil. Right: Experimental rig with an aerofoil in a wind tunnel.}
\label{fig:flutter_model}       % Give a unique label
\end{figure}

To provide a basis for our UDE model for aeroelastic flutter, we consider a physics-based dynamical model of an aerofoil (see the left panel of Fig.~\ref{fig:flutter_model}) \cite{Abdelkefi,Dimitriadis}. The wing, placed in the airstream with an airspeed of $V$, can move in the vertical direction and can rotate about the $\zeta$ axis. The motion of the aerofoil is characterised with the generalised coordinates $h$ (heave) and $\alpha$ (pitch angle), respectively. The wing is supported in the vertical direction with a linear spring and a damper, characterised by the stiffness $s$ and viscous damping $b$, while the compliance in the pitch is given by the nonlinear torsional stiffness $s_{\rm t}(\alpha) = s_{1 \rm t} \alpha + s_{3 \rm t} \alpha^3$ and linear torsional damping $b_{\rm t}$. This cubic stiffness term is a standard way to consider nonlinearity in low degrees-of-freedom aerofoil models \cite{Wright2015}. The bifurcation analysis of the model presented here is carried out in \cite{Abdelkefi}.

{To consider the effect of unsteady airflow around the aerofoil, two non-physical state variables: $w$ and its time-derivative $\dot w$ are introduced. This approximation of the effect of aerodynamic loads, assuming ideal, attached flow conditions, was first introduced by Theodorsen \cite{Theodorsen} and used nowadays as standard in low-order aeroelastic models \cite{Dimitriadis}.} %The result of the distributed aerodynamic load is given by the lift force $F$ and pitching moment $M$.

With the above considerations, in the deterministic case, the governing equations can be given in the form

\begin{equation}\label{eq:flutter_unsteady}
\mathbf{M}\ddot{\mathbf{y}} + \mathbf{B}\dot{\mathbf{y}} + \mathbf{S}{\mathbf{y}} +\mathbf{f}_{NL} (\mathbf{y}) = \mathbf{0},
\end{equation}
where the vector of generalised coordinates can be expressed as
\begin{equation}
\mathbf{y} = \left( \begin{matrix}
h \\
\alpha \\
w
\end{matrix} \right)
\end{equation}
whereas the vector of the nonlinear forces is given by
\begin{equation}
\mathbf{f}_{NL} (\mathbf{y}) = \left( \begin{matrix}
0 \\
s_{t2} \alpha^2 + s_{t3} \alpha^3 \\
0
\end{matrix} \right).
\end{equation}
The mass, damping and stiffness matrices $\mathbf{M}$, $\mathbf{B}$ and $\mathbf{S}$, containing both structural and aerodynamic terms are given in{~\ref{Flutter3}.}

In the following discussion, we will consider both an experimental and a numerical representation of the system as our observed system (the experimental version is practically relevant, but the numerical, while artificial, is helpful in assessing the merits of the proposed techniques). In the case of the latter, we consider the system given by the governing equations \eqref{eq:flutter_unsteady} as the \emph{observed} system. It is natural to assume that in a physical experiment the heave and pitch motion of the aerofoil would be relatively easy to measure, while the aerodynamic effects would remain hidden. Thus, the generalised coordinates $h$ and $\alpha$ describing the heave and pitch motions of the aerofoil are assumed to be \emph{measured}, while the non-physical generalised coordinate $w$ is treated as {\emph{unobserved}}.

Based on these assumptions, the UDE model of the system is formulated as

\begin{equation}\label{eq:flutter_uode}
\mathbf{M_{red}}\ddot{\mathbf{x}} + \mathbf{B_{red}}\dot{\mathbf{x}} + \mathbf{S_{red}}{\mathbf{x}} +\mathbf{f}_{NLred} (\mathbf{x}) + \mathbf{U}(\mathbf{x};V) = \mathbf{0},
\end{equation}
with
\begin{equation}
\mathbf{x} = \left( \begin{matrix}
h \\
\alpha
\end{matrix} \right),
\end{equation}
and
\begin{equation}
\mathbf{f}_{NLred} (\mathbf{y}) = \left( \begin{matrix}
0 \\
s_{t2} \alpha^2 + s_{t3} \alpha^3
\end{matrix} \right).
\end{equation}
The reduced mass, damping and stiffness matrices $\mathbf{M}_{\rm red}$, $\mathbf{B}_{\rm red}$ and $\mathbf{S}_{\rm red}$ are given in~\ref{Flutter2}. {Note that this model is the quasi-static approximation of the model given in Eq.~\eqref{eq:flutter_unsteady} assuming unsteady airflow. Nevertheless, while the two systems differ in dimension, both models can replicate bistable behaviour and loss of stability by a subcritical Hopf bifurcation in a qualitatively similar structure. Thus, the quantitative error between the two models is relatively small in this example.}

The universal approximator $\mathbf{U}$ that is considered to depend on the heave and pitch motion and the airspeed is used to compensate the error of the 2 DoF physics-based model compared to the 2+1 DoF model of the aerofoil.

To consider process noise in the system, Eq. \eqref{eq:flutter_unsteady} is extended by an additive noise term and reformulated as
\begin{equation}
\begin{split}
\rm{d} & \dot{\mathbf{z}} = - \mathbf{M}^{-1} \left( \mathbf{Bz} + \mathbf{Sy} + \mathbf{f}_{\rm{NL}}(\mathbf{y})\right) \rm{d}t - \mathbf{M}^{-1} \mathbf{F}_{\rm {L}} \rm{d}W, \\
\rm{d} & \mathbf{y} = \mathbf{z} \rm{d}t,
\end{split}
\end{equation}
where $\mathbf{F}_{\rm {L}}$ contains a moment acting on the aerofoil with a standard deviation $\sigma$
\begin{equation}
\mathbf{F}_{\rm {L}} = \left( \begin{matrix}
  0 \\ \sigma \\ 0
\end{matrix}
\right) .
\end{equation}
Note that, one would encounter process noise both in heave and pitch simultaneously in a real experiment. Thus, considering multiple sources of noise could result in a better characterisation of the system. However, we consider random perturbations in the pitch moment only, as this allows for a clearer quantification of the degree of noise-load in the presented examples.

Similarly, the corresponding UDE model is also considered in the SDE form
\begin{equation}\label{eq:flutter_sde}
  \begin{split}
\rm{d} & \dot{\mathbf{u}} = - \mathbf{M}_{\rm{red}}^{-1} \left( \mathbf{B}_{\rm{red}} \mathbf{u} + \mathbf{S}_{\rm{red}} \mathbf{x} + \mathbf{f}_{\rm{NLred}}(\mathbf{x}) + \mathbf{U}(\mathbf{x} ; V) \right) \rm{d}t - \mathbf{M}_{\rm{red}}^{-1} \tilde{\mathbf{F}}_{\rm {L}} \rm{d}W , \\
  \rm{d} & \mathbf{x} = \mathbf{u} \rm{d}t.
  \end{split}
\end{equation}
In view of the relatively small/moderately-scaled datasets used for training the machine-learnt structures, two types are considered: feed-forward, fully connected, relatively small-scale neural networks and Gaussian Processes.

\subsection{Nonlinear electrodynamic oscillator}\label{sec:energy_harv_num}
\begin{figure}
  \begin{center}
    \includegraphics[width=0.9\textwidth]{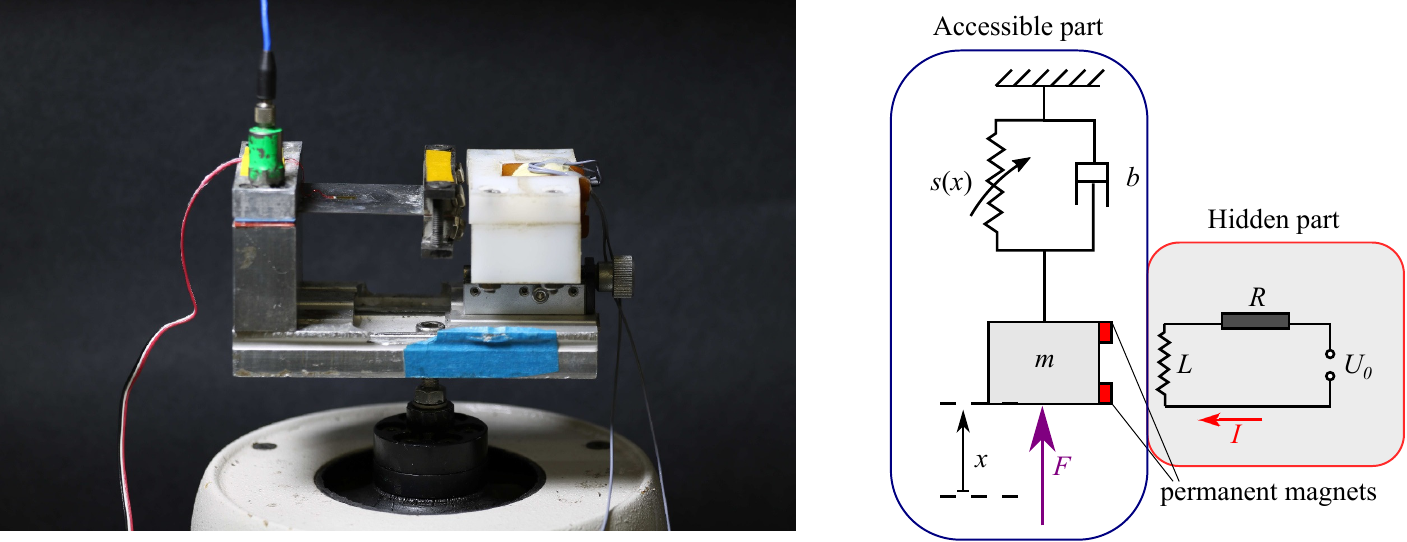}
  \end{center}
  \caption{The nonlinear electrodynamical oscillator (left panel) and its physics-based model (right panel). }
 \label{fig:eh_rig}
  \end{figure}

%\begin{figure}
%\caption{The nonlinear electrodynamical oscillator (left panel) and its physics-based model (right panel). }
%\label{fig:eh_rig}       % Give a unique label
%\end{figure}

\begin{figure}
  \begin{center}  
    \includegraphics[width=0.55\textwidth]{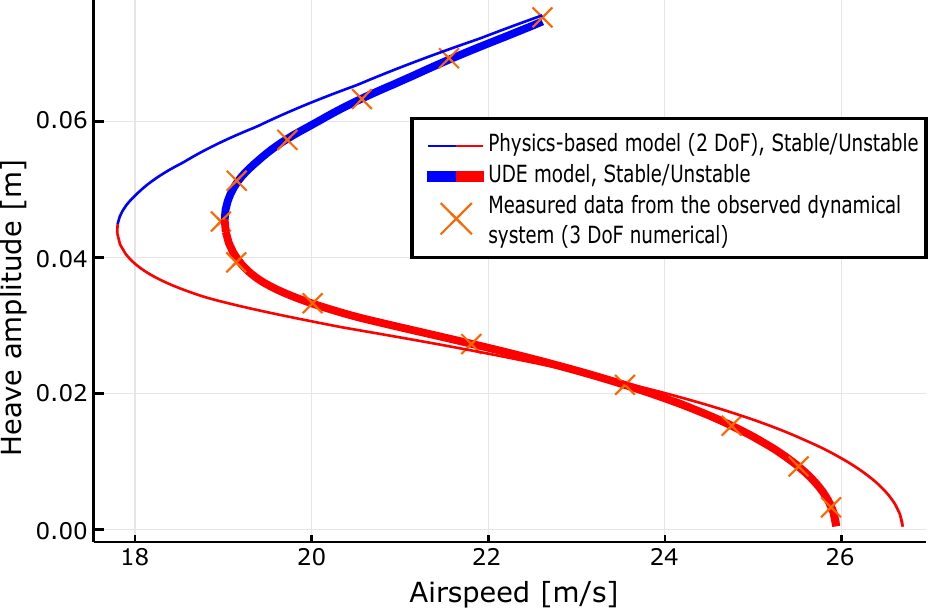}
  \end{center}
  \caption{Bifurcation diagrams of the UDE (thick curves) and pure mechanistic `steady-state' (thin curves) models of aeroelastic flutter against the training data generated by the deterministic, `unsteady' flutter model. In the UDE model, a neural network with 2 hidden layers of 12-12 neurons each, with a sigmoid activation function, was used to augment the physics-based part. The blue and red curves indicate stable and unstable limit cycles, respectively.}
  \label{fig:flutter_num_det}       % Give a unique label
  \end{figure}

% \begin{figure}
% % Use the relevant command to insert your figure file.
% % For example, with the graphicx package use
%   \includegraphics[width=0.45\textwidth]{duffing_model-eps-converted-to.pdf}
% % figure caption is below the figure
% \caption{Physics-based model of the nonlinear electrodynamic oscillator}
% \label{fig:eharvester}       % Give a unique label
% \end{figure}

Our second use case is a nonlinear electrodynamic oscillator subject to periodic forcing. In this case, we aim to identify a model for the physical experiment shown in Fig.~\ref{fig:eh_rig}. The nonlinear oscillator is formed from a thin steel plate, which is clamped to the base as a cantilever beam. At the end of the beam two iron masses, incorporating four permanent magnets, are attached. These magnets are interacting with an electromagnetic coil with an iron core in an insulated housing in the static part of the device. The resulting combination of structural and magnetic forces results in a nonlinear restoring force.

Based on a previous study on this particular device \cite{cbc_paper}, we consider an extended version of the Duffing oscillator with a seventh-order nonlinearity as the mechanistic model of the experiment (see Fig.~\ref{fig:eh_rig}). As such, we take the deformation of the beam and the corresponding speed as the measured state variables. The {undescribed} dynamics is mostly related to the electrodynamic coupling in the system. Our modelling framework however does not require a direct characterisation of the {undescribed} dynamics. It is worth mentioning that in \cite{Cammarano2011}, Cammarano et al. consider an LR circuit {(an electric circuit with a resistor and an inductor connected in series)} coupled to the mechanical part of the model which describes the system response to periodic excitation with good accuracy; however it did not accurately characterise the non-periodic response.

As shown in \cite{cbc_paper}, the model of the Duffing-type oscillator can capture the steady-state response to forcing amplitude at a constant frequency. However, it can only describe the frequency-response in a qualitative sense. Therefore, we consider a universal approximator depending on the forcing frequency $\omega$ alongside the physics-based part to generate our UDE model
\begin{equation}\label{eq:eh_ude}
  \dot y = - by - \alpha^2 x - c_3 x^3 - c_5 x^5 - c_7 x^7 + U(x,y;\omega) + A \cos (\omega t + \varphi),
\end{equation}
\begin{equation}
\dot x = y,
\end{equation}
where $x$ and $y$ are the displacement and speed of the mass $m$, $\alpha$ is the natural frequency of the oscillator, the parameter $b$ characterises the viscous damping, whereas the coefficients $c_3, c_5$ and $c_7$ describe the $3^{\rm rd}$-, $5^{\rm th}$- and $7^{\rm th}$-order nonlinearities, respectively. For practical reasons, the forcing amplitude $A$ is expressed with the amplitude $\Phi$ of the base acceleration as $A = c_A \omega_{\rm n}^2 \Phi$, where $\omega_{\rm n}$ is the linear natural angular frequency of the oscillator while $c_A$ is a constant coefficient. For this use-case, we study the deterministic case only; nevertheless, the model could also be expressed in the SDE form.

\section{UDE models of aeroelastic flutter trained {on a} numerical observed system}

\begin{figure*}
  \includegraphics[width=1.0\textwidth]{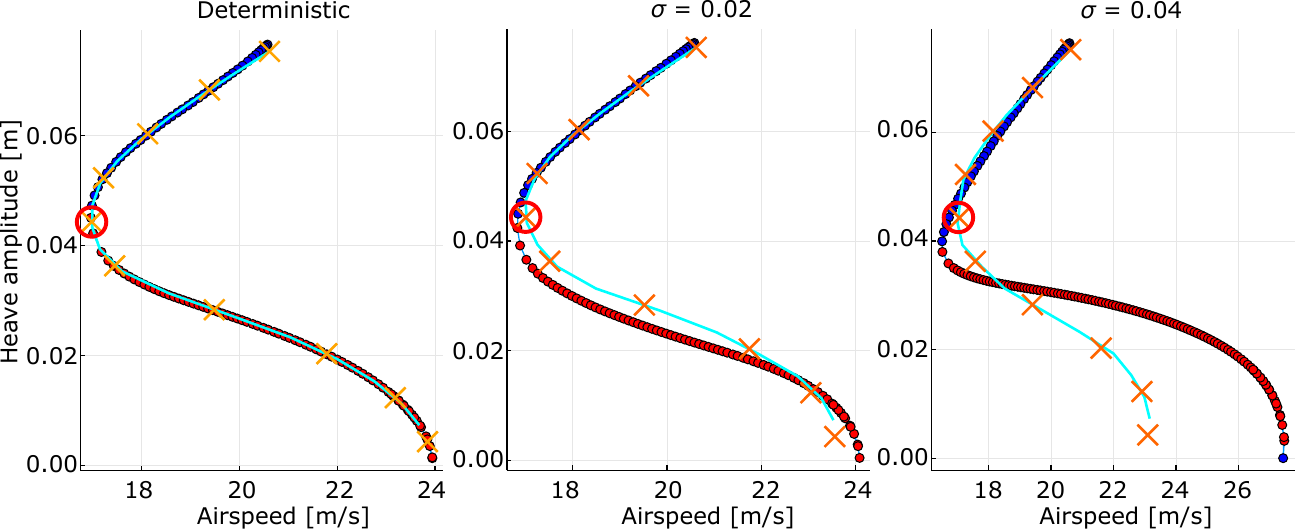}
\caption{Bifurcation diagrams of the pure mechanistic `steady-state' model of aeroelastic flutter and the UDE model using a Gaussian Process trained on data generated by the `unsteady' flutter model with process noise. The crosses indicate the average amplitudes of the near-periodic solutions of the stochastic underlying system while the predicted branches of limit cycles are generated using the deterministic parts of the UDE models.}
\label{fig:flutter_num_GP}       % Give a unique label
\end{figure*}

\begin{figure*}
  \includegraphics[width=1.0\textwidth]{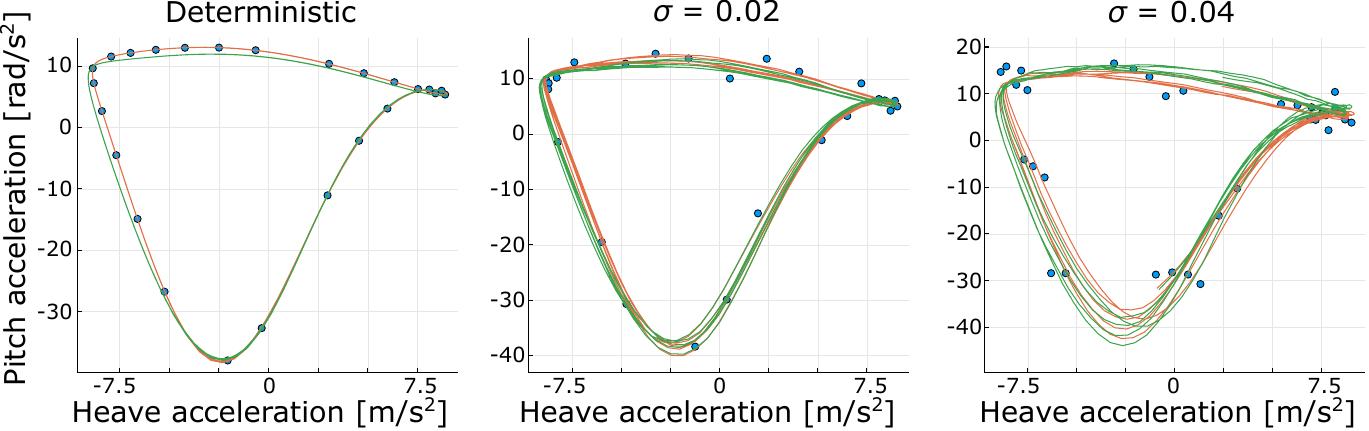}
\caption{Phase-portraits of the pitch and heave accelerations in the deterministic case and in the presence of process noise. The blue markers correspond to the `measured` data (underlying system) while the green and orange curves indicate the predictions with the pure mechanistic and the UDE models with a Gaussian Process, respectively. Each panel corresponds to the periodic or near-periodic solution marked with a red circle in the bifurcation diagrams in Fig.\ref{fig:flutter_num_GP}. }
\label{fig:flutter_pred_GP}       % Give a unique label
\end{figure*}

We carried out numerical simulations using the 2+1 DoF flutter model given in Eq.~\eqref{eq:flutter_unsteady}, to obtain its bifurcation diagram. We studied both the deterministic scenario and the case when the system was polluted with additive noise. The limit cycles (near-periodic solutions in the stochastic case) were traced using the technique of control-based continuation (CBC) \cite{SieberCBC,BartonCBC}. While in the deterministic case, a collocation-based method would be an alternative, the method of CBC can also handle stochastic systems, or even cases when the governing equations of the nonlinear system are not available, as is the case with physical experiments. The method is described briefly in {~\ref{App:CBC}}.

\subsection{Deterministic case}

The bifurcation diagrams of the different deterministic models of aeroelastic flutter are shown in Figure~\ref{fig:flutter_num_det}. The crosses indicate limit cycles of the unsteady physics-based model treated as the observed system. {Bifurcation diagrams may be used in practice as design aids to determine the airspeed-range where a particular aircraft is safe to fly. The model parameters were selected such that the system features a subcritical Hopf bifurcation generating a bistable airspeed range with a stable equilibrium co-existing with a stable limit cycle. This is not just the more challenging scenario in an experiment as, in contrast with the supercritical Hopf bifurcation, unstable limit cycles appear that are not trivial to find in experiments, but it is also the more dangerous case from the engineering point of view. In a real-life case, the bistable parameter-domain would be perceived as a sensitivity of the steady-state system behaviour to the initial conditions. Assuming that the stable equilibrium (i.e., no flutter) is the desired state of operation, it might be unsafe to operate an aircraft in the bistable speed-range as a sufficiently large perturbation could send the system into the domain of attraction of the large amplitude limit cycle, triggering flutter.}

As per the assumption of our modelling concept, we can identify a qualitatively correct reduced-order flutter model assuming a quasi-steady airflow using the measured state variables of the system (heave, pitch and the corresponding velocities) only. The bifurcation diagram of this model is indicated by the thin curve in Fig.~\ref{fig:flutter_num_det}. While this model is qualitatively correct in the sense that it has a bistable airspeed-range similar to the one of the observed system, one can also observe, that the predicted branch is slightly off the limit cycles captured in the `unsteady' flutter model. In contrast, the UDE model with a neural network with 2 hidden layers of 12-12 neurons each, with a sigmoid activation function, as indicated by the thick branch, is capable to capture the limit cycles almost perfectly.

As the difference between the physics-based model and the observed system is relatively small and we do not consider any random perturbation to the system either, even a fairly small-sized neural network provided satisfactory results. This model was also robust against the possible issues one may encounter during the training, i.e. we did not encounter the problem of overfitting while the identified minimum (albeit we do not have proof that it is the global one) always corresponded to a well-fitting model.

Arguably, this case with a relatively small quantitative error in the mechanistic model to compensate and no process noise is much simpler than a real experiment. Therefore, from a practical point of view it is also interesting to investigate the system with added process noise.

\subsection{Simulations with process noise}\label{sec:proc_noise}

\begin{figure}
\begin{center}
  \includegraphics[width=0.55\textwidth]{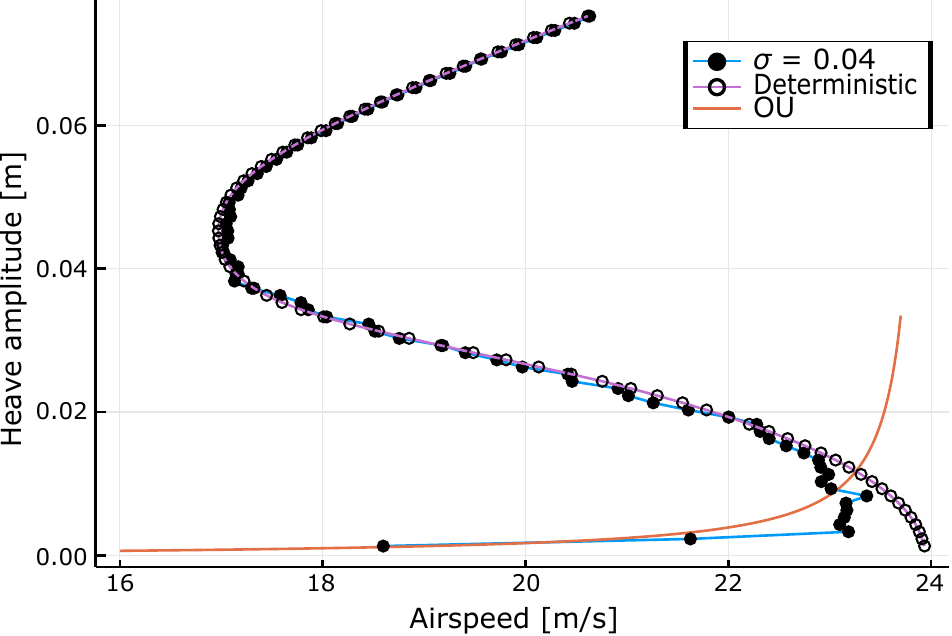}
\end{center}
\caption{ Branches of periodic/near-periodic solutions in the deterministic and stochastic flutter models as captured by the CBC algorithm. The orange curve indicates the vibration amplitude of the Ornstein-Uhlenbeck process corresponding to the `unsteady' flutter model linearised around the trivial equilibrium.}
\label{fig:OE}       % Give a unique label
\end{figure}

In the case where process noise is added to the numerical observed system, we employed Gaussian Processes as universal approximators alongside the physics-based part of the model.
The presence of process noise poses an additional challenge during the training procedure as the random perturbation results in near-periodic solutions instead of limit cycles, as would be the case in a deterministic system. In practice, this effect is perceived as a fluctuation in the vibration amplitude and frequency. The issue with using neural networks for stochastic systems exhibiting such behaviours is that, without any explicit condition to force the predicted trajectories to be periodic or explicitly identify the noise in the trajectories in the training data, there is no clear distinction whether the observed near-periodic trajectories are an inherent feature of the underlying deterministic system or purely a consequence of the stochastic perturbation. This might result in the deterministic neural network-based models yielding to near-periodic oscillations instead of limit cycles. Gaussian Processes are affected by this issue to a lesser degree as they handle the deterministic part of the model and the stochastic perturbation separately by definition.

Figure ~\ref{fig:flutter_num_GP} shows the branches of the identified limit cycles/near-periodic solutions with an additive stochastic moment with standard deviations $\sigma = 0, \quad 0.02$ and $0.04$ using the model given in \eqref{eq:flutter_sde}. The crosses corresponding to the `measurement data' indicate the average amplitude of the near-periodic solutions of the observed stochastic system, whereas the limit cycles of the UDE model were obtained by numerical collocation considering the deterministic part of the model only.

In the deterministic case ($\sigma = 0 $), the Gaussian process-based model can be fitted just as well to the data as it was shown with a neural network while it still provides a reasonably accurate model for a mildly perturbed system (see $\sigma = 0.02 $ ). Nevertheless, as the noise level increases, the fitted model becomes less accurate especially at lower vibration amplitudes where the signal to noise ratio is worse. From this point of view, it is a particular disadvantage of the Gaussian process-based models that they are fit directly to the acceleration-data rather than the trajectories and; as such, they are more severely affected by random perturbations. {This effect is illustrated by the phase plots in acceleration in Fig.~\ref{fig:flutter_pred_GP}}. This is the reason why the model accuracy can quickly deteriorate while only a mild perturbation can be observed in the heave and pitch time profiles.

While this problem might be alleviated by fitting the model based on the trajectories, as for the neural networks, using leave-one-out cross-validation \cite{Vehtari2016}, the Gaussian-Process-based models have another bottleneck as they do not scale well with larger amounts of data or more input variables. Even though the training procedure is still relatively fast on larger datasets, predictions are obtained slowly as the whole dataset is being used every time the right-hand-side of the UDE is evaluated. As a result, one may find it difficult to use the identified Gaussian-Process-based models for further analysis. Moreover, obtaining an accurate fit can be an issue as hyperparameters belonging to local minima may be more frequently found for more complex training data. Thus, we rather use neural network models to train a UDE model on experimental data. A possible alternative to this approach would be the use of sparse Gaussian Processes \cite{Titsias2009} which can eliminate the  issues caused by large datasets. The use of these structures however, is beyond the scope of this paper.

It is worth mentioning, that the presence of process noise results in a deviation from the deterministic solution-branch in the underlying model itself which is most prominent at lower vibration amplitudes. This can be best observed by comparing the critical airspeeds at which the Hopf bifurcation is detected for different noise levels in Fig.~\ref{fig:flutter_num_GP}. In the deterministic case, the Hopf bifurcation is approximately at 24 m/s whereas in the system where the pitch moment has a standard deviation of $\sigma = 0.04$ the critical airspeed is perceived around 23 m/s. This can be explained by fact that the noise-excitation leads to a near-periodic oscillation around the stable equilibria of the system \cite{SH_noise_stab}. This oscillation can be characterised by an Ornstein-Uhlenbeck process \cite{OE}, equivalent to the linearisation of the flutter model \eqref{eq:flutter_sde} around the equilibrium. In (see Fig.~\ref{fig:OE}) the near-periodic vibration amplitudes captured by CBC are compared to the Ornsten-Uhlenbeck model showing good agreement. Moreover, the perceived critical airspeed is very close to the point where the amplitude predicted by the Ornsten-Uhlenbeck model exceeds the amplitude of the unstable limit cycles of the deterministic system.

\section{Modelling an experimental observed system}

\subsection{Aeroelastic flutter}

% \begin{figure}
% % Use the relevant command to insert your figure file.
% % For example, with the graphicx package use
%   \includegraphics[width=0.45\textwidth]{Wind_tunnel_annotated.png}
% % figure caption is below the figure
% \caption{Experimental rig with an aerofoil in a wind tunnel.}
% \label{fig:flutter_exp_rig}       % Give a unique label
% \end{figure}

\begin{figure}
\begin{center}
  \includegraphics[width=0.95\textwidth]{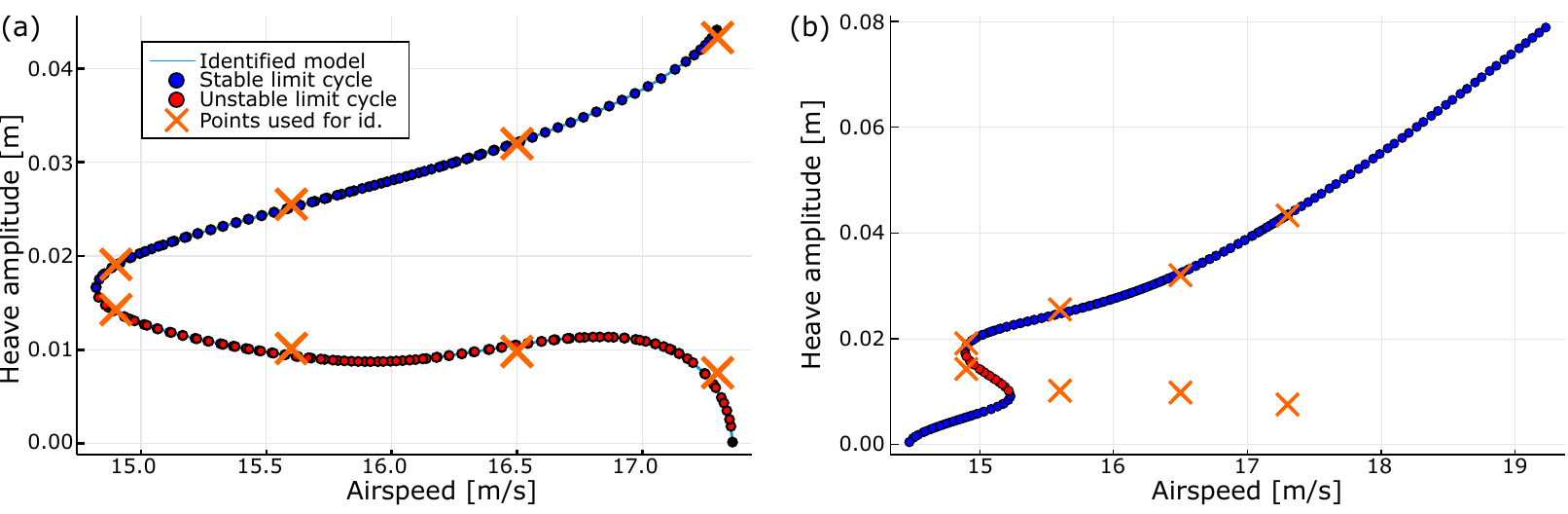}
\end{center}
\caption{Bifurcation diagrams of the UDE model of aeroelastic flutter using a neural network with 2 hidden layers of 48-48 neurons each with a sigmoid activation function, against the training data from physical experiment. The left panel shows a case when good fit to measurement data was achieved, whereas the right panel presents a bifurcation diagram of a model corresponding to a local minimum of the objective function.}
\label{fig:flutter_exp_bif}       % Give a unique label
\end{figure}

% \begin{figure}
% % Use the relevant command to insert your figure file.
% % For example, with the graphicx package use
%   \includegraphics[width=0.45\textwidth]{flutter_exp_bif_localmin-eps-converted-to.pdf}
% % figure caption is below the figure
% \caption{Bifurcation diagram of the UDE model using a neural network with 2 hidden layers of 48-48 neurons each with a sigmoid activation function, corresponding to a local minimum of the objective function against the training data from physical experiment.}
% \label{fig:flutter_exp_bif_localmin}       % Give a unique label
% \end{figure}

\begin{figure}[h!]
  \begin{center}
    \includegraphics[width=0.95\textwidth]{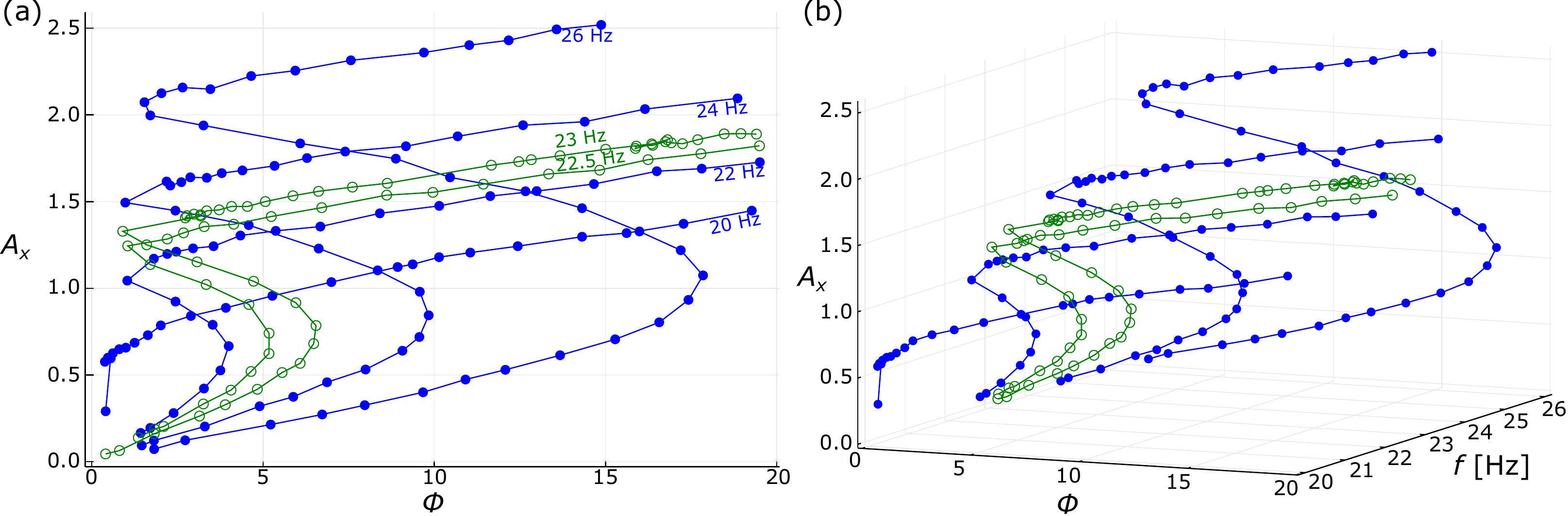}
  \end{center}
  \caption{The bifurcation diagrams used as training data, collected from the physical nonlinear oscillator (Panels (a) and (b) show the same data in 2D and 3D, respectively). The dataset with a narrow frequency range (22.5-23 Hz) is shown is green with empty markers, while the blue curves and filled markers correspond to the wide frequency range dataset (20-26 Hz)}
  \label{fig:eh_data_2D}       % Give a unique label
\end{figure}

\begin{figure*}
    \includegraphics[width=1.0\textwidth]{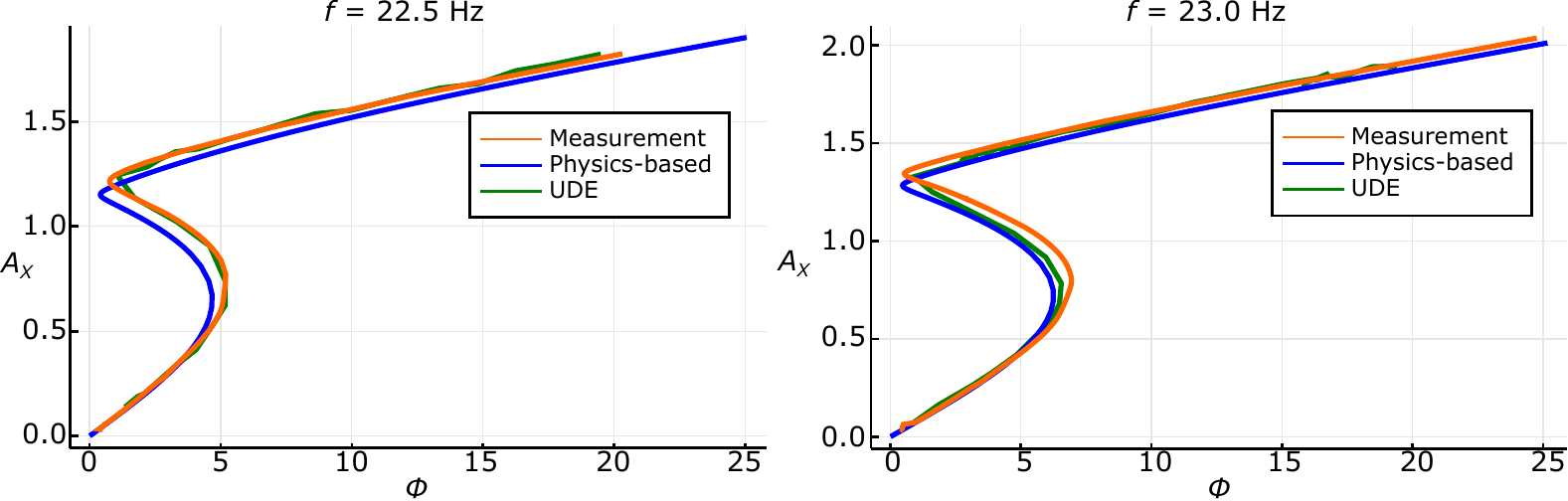}
  \caption{Bifurcation diagrams of the pure mechanistic and the UDE models against the measurement data (narrow frequency-range dataset). In the UDE model, the mechanistic part was augmented with a neural network of 2 hidden layers of 48-48 neurons each with a sigmoid activation function.}
  \label{fig:eh_narrow_bif}       % Give a unique label
\end{figure*}
    
\begin{figure*}
      \includegraphics[width=1.0\textwidth]{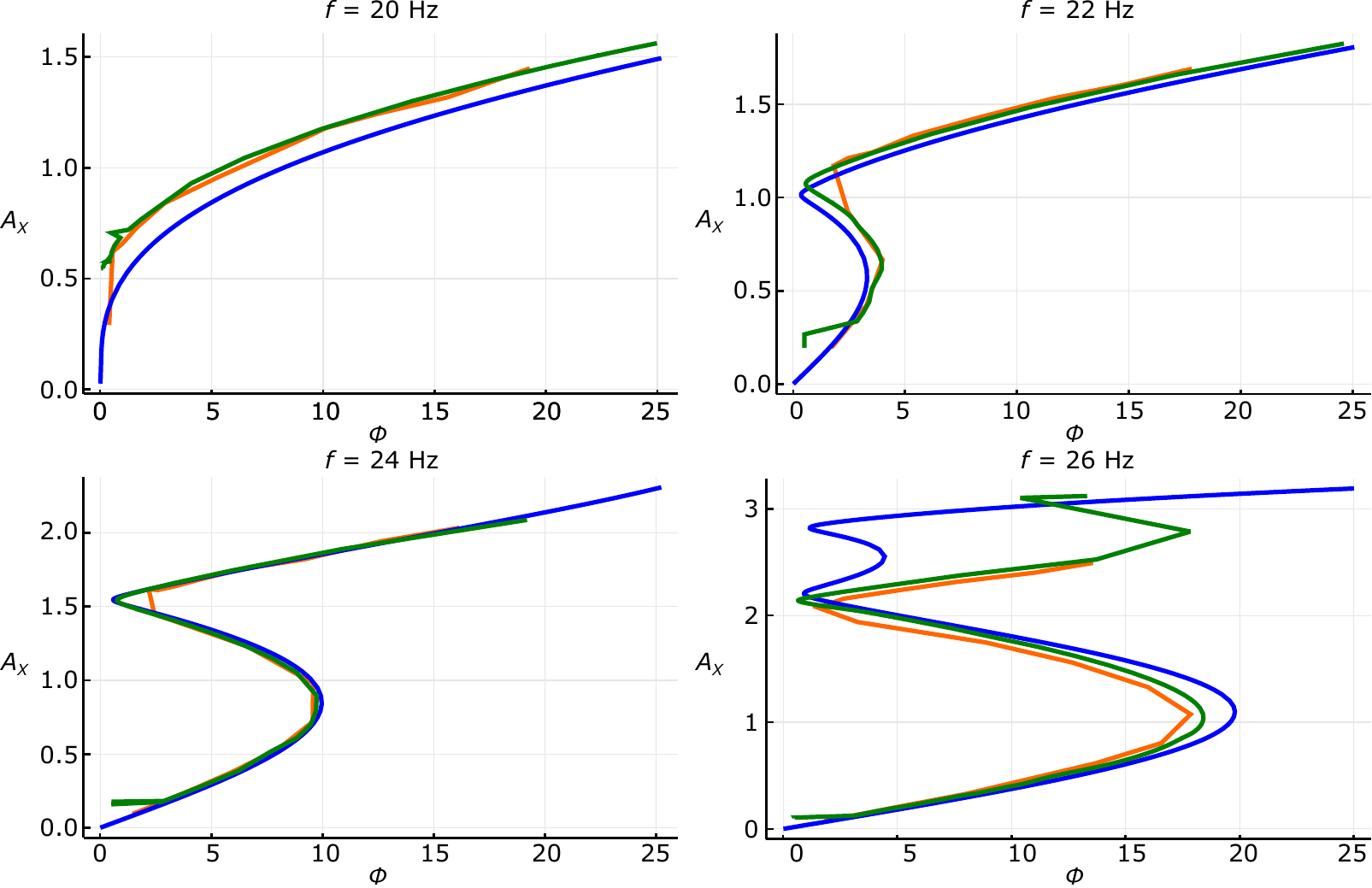}
    \caption{Bifurcation diagrams of the UDE and pure mechanistic models against the measurement data (wider frequency-range dataset). In the UDE model, the mechanistic part was augmented with a neural network of 2 hidden layers of 48-48 neurons each with a sigmoid activation function.}
    \label{fig:eh_wide_bif}       % Give a unique label
\end{figure*}

We also used the UDE model given in \eqref{eq:flutter_uode} of the aerofoil to model the dynamics of a physical experiment with a wing profile in the wind tunnel shown in the right panel of Fig.~\ref{fig:flutter_model}. The architecture of the model was the same as in case of the numerical simulations, i.e. pitch and heave (and their respective time derivatives) were assumed to be the measured state variables while we considered everything else as the hidden dynamics of the experiment. {This example relies on the measurement data published by Lee et al. \cite{Lee2020} who use it to find a reduced-order model of the experiment. The measurement dataset is available at \cite{Kyoungdata}. It is worth mentioning that the data from the wind tunnel test rig was collected as proof of concept for control-based continuation (CBC) in a challenging experiment where the vibration is self-excited and significant process noise is present. Thus, the experimental regime was not specifically tailored for training UDE models.}

This problem proved to be too complex to be modelled using Gaussian processes, as the size of the training dataset that adequately represents the observed dynamics exceeds the size where these structures are practical to use. As a result, predictions with Gaussian-Process-based models get very costly computationally making the evaluation and further use of the trained models infeasible. Thus, for the physical experiment, we employ neural networks as machine learnable structures.

Similarly to the numerical case, we use time series data from the limit cycles to train the neural networks. As explained in Section~\ref{sec:proc_noise}, it is a strong requirement for the training data to comprise limit cycles, as neural-network-based models are prone to replicate the near-periodic oscillations observed in a system with random excitation. Thus, the measurement data had to be pre-processed to generate quasi-deterministic, periodic time profiles from the otherwise noise-polluted measurement data. This was achieved by Fourier-transformation and averaging over longer data-segments at each measurement points. As a result, we generated 8 time-series of limit cycles corresponding to different airspeeds, 4 solutions taken from the stable and the unstable branches each.

The left panel of Fig.~\ref{fig:flutter_exp_bif} shows the bifurcation diagram of the trained UDE model which fits well to the measurement data. Unfortunately however, repeated attempts to train the same model starting with a different initial guess of the neural network's gains and biases on the same data, often resulted in poorly-fitting models as the training algorithm was prone to converge to local minima of the objective function instead of the global one (see e.g. the right panel of Fig.~\ref{fig:flutter_exp_bif}). The robustness of the training algorithm could be improved by techniques such as drop-out or mini-batching \cite{Shreshta2019NN} and presumably also by adding more data. Nevertheless, one may never fully eliminate the sensitivity to the initial guess and choice of algorithm parameters. This highlights one of the main conclusions of our study; namely, within the current framework, one should always check the adequacy of the result rather than just accepting a trained UDE model.

One of the main weakness of the current training framework using the DiffEqFlux.jl package is that by using forward simulations to generate the predictions we never use the knowledge that the measured trajectories, used as training data, correspond to limit cycles. While the error measure, i.e. the Euclidean distance between the observed and predicted trajectories is adequate in a sense that this should be minimal if the prediction is accurate, the approximation error may get amplified in certain cases, i.e. near to an unstable limit cycle. Ultimately, the success of the training procedure largely depends on the complexity of the problem and the ratio of the known and hidden dynamics. The training procedure for systems with limit cycles potentially could be made more robust if the training algorithm was incorporated into a boundary-value problem solver rather than a simulation of an initial value problem.

\subsection{Nonlinear electrodynamical oscillator}

% \begin{figure}
% % Use the relevant command to insert your figure file.
% % For example, with the graphicx package use
%   \includegraphics[width=0.5\textwidth]{exp_data_eharv-eps-converted-to.pdf}
% % figure caption is below the figure
% \caption{3-dimensional diagram of steady-state response of the physical nonlinear oscillator. The dataset with a narrow frequency range (22.5-23 Hz) is shown in green with empty markers, while the blue curves and filled markers correspond to the wide frequency range dataset (20-26 Hz)}
% \label{fig:eh_data}       % Give a unique label
% \end{figure}

The UDE model of the electrodynamic oscillator given in Eq.~\eqref{eq:eh_ude} is trained on measurement data collected from the experimental rig shown in Fig.~\ref{fig:eh_rig}. The steady-state limit cycles of the observed system were traced at several forcing frequencies while the forcing amplitude was varied. Thus, by the continuation of the S-shaped response curves, corresponding to the individual frequencies, we obtained the response-surface shown in Fig.~\ref{fig:eh_data_2D}. {In \cite{cbc_paper}, a similar dataset, comprising stable and unstable limit cycles in the system, is used to assess the robustness of CBC to noise by identifying the parameters of a Duffing-type oscillator model.}

{The structure of the solutions indicate a bistable behaviour with a low and a large amplitude stable limit cycle and an unstable limit cycle in the system}.
To assess the performance of the hybrid model, we consider two datasets, one a with a narrow (22.5-23 Hz) and another with a wide frequency range (20-26 Hz with an increment of 2 Hz) indicated by the green and blue curves, respectively.

Figures \ref{fig:eh_narrow_bif} and \ref{fig:eh_wide_bif} show the measured bifurcation diagrams against the predictions with the UDE model and the bare physics-based model (with the parameters identified using the system response at 24 Hz) for the narrow and wide-band cases respectively. One can observe that the addition of the neural network to the mechanistic model clearly results in a more accurate prediction; however, the extent of the improvement is not consistent for the two datasets.

In case of the narrow-frequency-range dataset, the UDE model fits very well to the measurement data at 22.5 Hz and provides a considerably better prediction for the amplitudes at 23 Hz too, even though there is still some error in the estimation.

Notably, we have not found an UDE model showing similar performance on the wide-frequency-range dataset as the extent of the remaining error between the measured bifurcation diagrams and the predictions is somewhat larger than for the narrow-frequency-range dataset. Moreover, in this case, in each effort to train the UDE model, we encountered the problem of overfitting; i.e. the resulting hypersurface of the underlying function of the neural network had large variations. This effect made the identified model difficult to handle. Also, we were unable to trace the limit cycles with a conventional collocation-based method. We could overcome this issue by applying the CBC technique, which is more robust against the random variations, on the trained UDE model. The result of using CBC indicates that in spite of the observed overfitting, the UDE model still makes an improvement in the accuracy of the prediction. Even though the wide-frequency-range problem seems to be slightly beyond the limits in terms of complexity our approach can handle within the current training framework, the narrow-bandwith case demonstrates the potential of the concept.

%It is worth to assess a popular way to regularise neural networks \cite{Shreshta2019NN} prone to overfitting by adding the norm of the weights and biases to the objective function. Although this `trick' results in a clearly smoother underlying function of the neural network one also compromises the accuracy of the model (see Fig...). Thus, this cannot be seen as a general solution to this problem.

\section{Conclusions}

Our study exploited examples of using universal approximators alongside mechanistic models of engineering systems with the aim to augment the qualitatively correct models to reach quantitative accuracy. Our trials demonstrate the potential and some challenges with the technique, underlining that the complexity of the problem and the quality of the training data is the most important factor in the training procedure. This by itself highlights the relevance of using hybrid models instead of differential equations with machine-learnable structures only (without a physics-based part), since in general, the closer the qualitative models are to the `ground truth' the easier it is to capture the `missing part' of the dynamics with machine learnable structures.

Nevertheless, in some cases we identified issues like converging to local minima or overfitting during the training procedure. While it is not always possible to identify the source of these issues, one can generally say that problem complexity or insufficient training data plays an important part in it. Therefore, further efforts need to be invested in finding what type of training data could enhance the robustness of the training procedure. Due to these issues, within the current framework, one should be critical to the trained UDE models and an a-posteriori assessment is needed to establish whether they are adequately represent the modelled system. In other words, the training process is not robust enough for its results to be blindly accepted.

Regarding the performance of the two different types of universal approximators used in our analysis, we can generally conclude that while Gaussian processes can perform better in systems with a mild noise load, neural networks were more versatile, as they could handle more complex problems as well.

Our study mainly focuses on exploring how the most common approaches to train UDE models perform in replicating bifurcation diagrams of nonlinear dynamical systems and identifying the challenges that may arise; yet, one has to also admit that, arguably, there is still potential in developing the training methods to optimise their performance for nonlinear dynamical systems with limit cycles. From this point of view, trying to match the right-hand side of the UDEs to the acceleration data, as done in the examples where Gaussian processes were used as the universal approximator, is a somewhat naive approach as this does not reflect at all to the fact that dynamical systems are considered. From this point of view, the simulation-based approach of the DiffEqFlux.jl package in Julia, we deployed for the neural-network-based models, is better. Yet, not even this method uses the periodicity of the trajectories used as training data in a direct way. Instead, predictions are provided as solutions of initial value problems assuming that if the predictions are accurate the resulting model will also have a limit cycle. In \cite{Roesch2021}, it was demonstrated that collocation can make the training algorithms more robust in retaining the qualitative features of the observed dynamical system. Thus, implementing machine learning techniques using boundary-value problem solvers, better suited to find periodic solutions than time-forward simulations, could potentially improve the robustness of the fitting procedure in a significant way.

Nevertheless, our exploratory study has shown that using mechanistic-machine learnable hybrid models for varying-parameter problems of nonlinear systems has clearly a potential to deliver accurate, low-dimensional models of these systems.

\section*{Acknowledgements}
This research has received funding from the {\em Digital twins for improved dynamic design} (EP/R006768/1) EPSRC grant. The support of the EPSRC is greatly acknowledged.
%\end{acknowledgement}


%\bibliographystyle{spphys}
\bibliographystyle{unsrt}      % mathematics and physical sciences
\bibliography{ML_hybrid.bib}   % name your BibTeX data base

%\begin{thebibliography}{}
% \bibitem{DiffEqFlux}
% % Format for Journal Reference
% Rackauckas C., Ma Y., Martensen J., Warner C., Zubov K., Supekar R., Skinner D., Ramadhan A., Edelman A., Universal Differential Equations for Scientific Machine Learning, preprint (2020), https://arxiv.org/abs/2001.04385
% Author, Book title, page numbers. Publisher, place (year)
% \end{thebibliography}

\appendix

\section{Coefficient matrices of the flutter model}

\setcounter{equation}{0}
\renewcommand{\theequation}{A.\arabic{equation}}

\subsection{Unsteady model}\label{Flutter3}

The mass matrix $\mathbf{M}$ in Eq.~\eqref{eq:flutter_unsteady} is divided into a structural and an aerodynamic part
\begin{equation}
  \mathbf{M} = \mathbf{M}_{\rm st} + \mathbf{M}_{\rm ae},
\end{equation}
where
\begin{equation}
  \mathbf{M}_{\rm st} = \left( \begin{matrix}
    m & m_{\rm w}(b-x_f) & 0 \\
    m_{\rm w}(b-x_f) & I_{\rm A} & 0 \\
    0 & 0 & 0 \\
  \end{matrix} \right),
\end{equation}
and
\begin{equation}
  \mathbf{M}_{\rm ae} = \left( \begin{matrix}
    \rho b^2 \pi & -\rho b^3 \pi a & 0 \\
    -\rho b^3 \pi a & \rho b^4 \pi (\frac{1}{8}+a^2) & 0 \\
    0 & 0 & 1 \\
  \end{matrix} \right)
\end{equation}
Similarly, the damping and stiffness matrices $\mathbf{B}$ and $\mathbf{S}$ are given in the same form
\begin{equation}
  \mathbf{B} = \mathbf{B}_{\rm st} + \mathbf{B}_{\rm ae},
\end{equation}
where
\begin{equation}
  \mathbf{B}_{\rm st} = \left( \begin{matrix}
    d & 0 & 0 \\
    0 & d_{\rm t} & 0 \\
    0 & 0 & 0 \\
  \end{matrix} \right), \quad  \mathbf{B}_{\rm ae} = \left( \begin{matrix}
    b_{\rm{ae11}} & b_{\rm{ae12}} & b_{\rm{ae13}} \\
    b_{\rm{ae21}} & b_{\rm{ae22}} & b_{\rm{ae23}} \\
    -\frac{1}{b} & -(\frac{1}{2}-a) & (c_2+c_4) \frac{V}{b} \\
  \end{matrix} \right),
\end{equation}
% \begin{equation}
%   \mathbf{B}_{\rm ae} = \left( \begin{matrix}
%     b_{\rm{ae11}} & b_{\rm{ae12}} & b_{\rm{ae13}} \\
%     b_{\rm{ae21}} & b_{\rm{ae22}} & b_{\rm{ae23}} \\
%     -\frac{1}{b} & -(\frac{1}{2}-a) & (c_2+c_4) \frac{V}{b} \\
%   \end{matrix} \right),
% \end{equation}
and
\begin{small}
\begin{multline}
    b_{\rm{ae11}} = 2 \pi \rho b V (c_0-c_1-c_3), \; b_{\rm{ae12}} = \rho b^2 \pi V+2 \pi \rho b^2 V (\frac{1}{2}-a) (c_0-c_1-c_3),\\
    b_{\rm{ae13}} = 2 \pi \rho V^2 b (c_1 c_2+c_3 c_4), \; b_{\rm{ae21}} = -2 \pi \rho b^2 V (a+\frac{1}{2}) (c_0-c_1-c_3), \\
    b_{\rm{ae22}} = \rho b^3 \pi V (\frac{1}{2}-a)-2 \rho b^3 V \pi (\frac{1}{2}-a) (\frac{1}{2}+a) (c_0-c_1-c_3), \\
    b_{\rm{ae23}} = -2 \pi \rho b^2 V^2 (a+\frac{1}{2}) (c_1 c_2+c_3 c_4). \\
  \end{multline}
\end{small}
% \begin{equation}
% b_{\rm{ae11}} = 2 \pi \rho b V (c_0-c_1-c_3),
% \end{equation}
% \begin{equation}
% b_{\rm{ae12}} = \rho b^2 \pi V+2 \pi \rho b^2 V (\frac{1}{2}-a) (c_0-c_1-c_3),
% \end{equation}
% \begin{equation}
% b_{\rm{ae13}} = 2 \pi \rho V^2 b (c_1 c_2+c_3 c_4),
% \end{equation}
% \begin{equation}
% b_{\rm{ae21}} = -2 \pi \rho b^2 V (a+\frac{1}{2}) (c_0-c_1-c_3),
% \end{equation}
% \begin{equation}
% b_{\rm{ae22}} = \rho b^3 \pi V (\frac{1}{2}-a)-2 \rho b^3 V \pi (\frac{1}{2}-a) (\frac{1}{2}+a) (c_0-c_1-c_3),
% \end{equation}
% \begin{equation}
% b_{\rm{ae23}} = -2 \pi \rho b^2 V^2 (a+\frac{1}{2}) (c_1 c_2+c_3 c_4).
% \end{equation}
The stiffness matrix can be expressed as
\begin{equation}
  \mathbf{S} = \mathbf{S}_{\rm st} + \mathbf{S}_{\rm ae},
\end{equation}
with
\begin{equation}
  \mathbf{S}_{\rm st} = \left( \begin{matrix}
    s & 0 & 0 \\
    0 & s_{\rm t} & 0 \\
    0 & 0 & 0 \\
  \end{matrix} \right), \quad  \mathbf{S}_{\rm ae} = \left( \begin{matrix}
    0 & s_{\rm{ae12}} & s_{\rm{ae13}} \\
    0 & s_{\rm{ae22}} & s_{\rm{ae23}} \\
    0 & -\frac{V}{b} & c_2 c_4 \frac{V^2}{b^2}\\
  \end{matrix} \right),
\end{equation}
% and
% \begin{equation}
%   \mathbf{S}_{\rm ae} = \left( \begin{matrix}
%     0 & s_{\rm{ae12}} & s_{\rm{ae13}} \\
%     0 & s_{\rm{ae22}} & s_{\rm{ae23}} \\
%     0 & -\frac{V}{b} & c_2 c_4 \frac{V^2}{b^2}\\
%   \end{matrix} \right),
% \end{equation}
where
\begin{multline}
    s_{\rm{ae12}} = 2 \pi \rho b V^2 (c_0-c_1-c_3), \; s_{\rm{ae13}} = 2 \pi \rho V^3 c_2 c_4 (c_1+c_3), \\
    s_{\rm{ae22}} = -2 \pi \rho b^2 V^2 (\frac{1}{2}+a) (c_0-c_1-c_3), \\ s_{\rm{ae23}} = -2 \pi \rho b^2 V^3 (\frac{1}{2}+a) c_2 c_4 (c_3+c_1). \\
\end{multline}
% \begin{equation}
% s_{\rm{ae12}} = 2 \pi \rho b V^2 (c_0-c_1-c_3),
% \end{equation}
% \begin{equation}
% s_{\rm{ae13}} = 2 \pi \rho V^3 c_2 c_4 (c_1+c_3),
% \end{equation}
% \begin{equation}
% s_{\rm{ae22}} = -2 \pi \rho b^2 V^2 (\frac{1}{2}+a) (c_0-c_1-c_3),
% \end{equation}
% \begin{equation}
% s_{\rm{ae23}} = -2 \pi \rho b^2 V^3 (\frac{1}{2}+a) c_2 c_4 (c_3+c_1).
% \end{equation}

\subsection{Quasi-steady airflow model}\label{Flutter2}
The coefficient matrices in the equation of motion of the 2 DoF flutter model (see Eq. \eqref{eq:flutter_uode}) can be expressed as follows. The mass matrix $\mathbf{M}_{\rm red}$ reads
\begin{equation}
  \mathbf{M}_{\rm red} = \mathbf{M}_{\rm red, st} + \mathbf{M}_{\rm red, ae},
\end{equation}
where
\begin{equation}
  \mathbf{M}_{\rm red, st} = \left( \begin{matrix}
    m & m_{\rm w}(b-x_f)\\
    m_{\rm w}(b-x_f) & I_{\rm A}\
  \end{matrix} \right),
\end{equation}
and
\begin{equation}
  \mathbf{M}_{\rm red, ae} = \left( \begin{matrix}
    \rho b^2 \pi & -\rho b^3 \pi a\\
    -\rho b^3 \pi a & \rho b^4 \pi (\frac{1}{8}+a^2)\\
    \end{matrix} \right).
\end{equation}
The damping matrix is given by
\begin{equation}
  \mathbf{B}_{\rm red} = \mathbf{B}_{\rm red, st} + \mathbf{B}_{\rm red, ae},
\end{equation}
\begin{equation}
  \mathbf{B}_{\rm red, st} = \left( \begin{matrix}
    d & 0\\
    0 & d_{\rm t}\\
  \end{matrix} \right), \; \mathbf{B}_{\rm red, ae} = \left( \begin{matrix}
    2 \pi \rho b V & 2 \pi \rho b^2 V (1-a) \\
    -2 \pi \rho b^2 V (a+\frac{1}{2}) & -2 a \rho b^3 V \pi (\frac{1}{2}-a) \\
  \end{matrix} \right),
\end{equation}
and
% \begin{equation}
%   \mathbf{B}_{\rm red, ae} = \left( \begin{matrix}
%     2 \pi \rho b V & 2 \pi \rho b^2 V (1-a) \\
%     -2 \pi \rho b^2 V (a+\frac{1}{2}) & -2 a \rho b^3 V \pi (\frac{1}{2}-a) \\
%   \end{matrix} \right),
% \end{equation}
whereas the stiffness matrix reads
\begin{equation}
  \mathbf{S}_{\rm red} = \mathbf{S}_{\rm red, st} + \mathbf{S}_{\rm red, ae},
\end{equation}
with
\begin{equation}
  \mathbf{S}_{\rm red, st} = \left( \begin{matrix}
    s & 0\\
    0 & s_{\rm t}\\
  \end{matrix} \right), \; \mathbf{S}_{\rm red, ae} = \left( \begin{matrix}
    0 & 2 \pi \rho b V^2 \\
    0 & -2 \pi \rho b^2 V^2 (\frac{1}{2}+a) \\
  \end{matrix} \right).
\end{equation}
and
% \begin{equation}
%   \mathbf{S}_{\rm red, ae} = \left( \begin{matrix}
%     0 & 2 \pi \rho b V^2 \\
%     0 & -2 \pi \rho b^2 V^2 (\frac{1}{2}+a) \\
%   \end{matrix} \right).
% \end{equation}

\section{Tracing periodic solutions in the flutter model with control-based continuation} \label{App:CBC}

\setcounter{equation}{0}
\renewcommand{\theequation}{B.\arabic{equation}}

The limit cycles, we use as training data, are tracked with the technique of control-based continuation (CBC) for both the physical and numerical simulations. As shown in previous studies \cite{cbc_paper}, this method is fairly robust against random perturbations. Thus, it can be used to survey the underlying deterministic system in noise-polluted experiments.

CBC is based on applying a stabilising and non-invasive control to the experimental system, such that the system converges to a periodic solution and, by iterating the control target, the algorithm finds a limit cycle that is also a solution of the open/loop (uncontrolled) system.

In case of the energy harvester experiments the open-loop system is subject to harmonic forcing, whereas the fluttering aerofoil is an autonomous system. Hence, the period of the limit cycles in the flutter case is an additional unknown parameter. The algorithm used to track the limit cycles of the physical nonlinear oscillator is described in detail in \cite{cbc_paper}. In this section, we present the algorithm used for the numerical simulation of the fluttering aerofoil. This is an extended version of the algorithm presented in \cite{Lee2020} by Lee et al., capable of finding limit cycles in autonomous system. However, due to the limitations of the experimental rig, limit cycles are only tracked for a number of fixed airspeeds while we implement a full continuation using the airspeed as bifurcation parameter.

To stabilise the system at a steady-state limit cycle, a control force is applied to the aerofoil
\begin{equation}\label{eq:cbc_ctrl}
F_{\rm{ctrl}} = k_{\rm p} (h_{\rm{t}}(t) - h(t)) + k_{\rm d} (\dot h_{\rm{t}}(t) - \dot h(t)),
\end{equation}
where $h_{\rm{t}}(t)$ is the heave control target whereas $k_{\rm p}$ and $k_{\rm d}$ are the proportional and derivative control gains. Thus, in the deterministic case, the equation of motion \eqref{eq:flutter_unsteady} is extended by a forcing term in the right-hand-side
\begin{equation}
\mathbf{M}\ddot{\mathbf{y}} + \mathbf{B}\dot{\mathbf{y}} + \mathbf{S}{\mathbf{y}} +\mathbf{f}_{NL} (\mathbf{y}) = \mathbf{F},
\end{equation}
with $\mathbf{F} = ( \begin{matrix} F_{\rm{ctrl}} & 0 & 0 \end{matrix})^{\rm T}$. To consider process noise, we apply the same extension to the SDE in \eqref{eq:flutter_sde}.

To consider a control target that is phase-locked to the self-excited vibration of the system, the instantaneous phase, defined by
\begin{equation}
\phi (t) = {\rm arctan} \left( \dot h (t), \omega h (t) \right),
\end{equation}
is used with $\omega$ denoting the measured angular frequency. Then, the control target is considered as
\begin{equation}
h_{\rm t} (t) = A_{\rm 0t} + A_{\rm 1t} \cos \left( \phi(t) \right),
\end{equation}
where $A_{\rm 0t}$ and $A_{\rm 1t}$ are the target static deflection and amplitude, respectively. The time derivative of the instantaneous phase is approximated with the help of the measured angular frequency as $\dot \phi \approx -\omega$; thus, the time derivative of the control target is given by
\begin{equation}
\dot h_{\rm t} (t) = \omega A_{\rm 1t} \sin \left( \phi(t) \right).
\end{equation}

\subsection{Root-finding and continuation}

If the control law in Eq.~\eqref{eq:cbc_ctrl} is stabilising, the system settles at a steady state periodic solution in the deterministic case, or a quasi-periodic solution in the presence of process noise. This solution is considered in a similar form to the control target
\begin{equation}
h (t) \approx A_{0} + A_{1} \cos \phi(t).
\end{equation}
In case of the near-periodic solutions in the stochastic system, we average over several periods, and the solution given by average coefficients are accepted as the steady-state periodic solution. If we find the control target coefficients $A_{\rm 0t}^\ast$ and $A_{\rm 1t}^{\ast}$ where the steady-state errors $e_0 = A_{\rm 0t} - A_0$, $e_1 = A_{\rm 1t} - A_1$ yield to zero (or are below a defined tolerance), this solution is accepted as a limit cycle in the open-loop system. While testing the algorithm, we found that the fixed point iteration $A_{\rm 0t}^{k+1}:=A_{0}^{k}$ was enough to eliminate the error $e_0$ in the static deflection; hence, we only consider the error in the vibration amplitude $e_1$ in continuation algorithm.

For the continuation of the branch of limit cycles, we consider the root-finding problem in two variables, the target amplitude and the airspeed $V$
\begin{equation}\label{eq:err1}
  e_1(A_{\rm 1t}, V) = A_{\rm 1t} - A_1(V).
\end{equation}

We take advantage from the fact, that a unique limit cycle corresponds to every vibration amplitude in the flutter model. Therefore, we assume that Eq. \eqref{eq:err1} has a unique solution $V^{\ast}$ for all target vibration amplitudes. As in our case, the stochastic nature of the system may prohibit the implementation of a direct derivative-based (e.g. Newton-like) root-finding algorithm for this problem, we make a mesh in airspeeds $V$ such that it includes the suspected solution $V_{\rm{min}}=V_0<V_1<...<V_n=V_{\rm{max}}$. Then, we register the corresponding steady-state vibration amplitudes $A_1^0$, $A_1^1$, ... $A_1^n$ and fit a $3^{\rm{rd}}$-order polynomial to these samples and root-finding is performed on this surrogate model. After identifying a limit cycle, the continuation progresses by incrementing/decrementing the target vibration amplitude $A_{\rm 1t}$.

\end{document}